\documentclass{article}
\usepackage[utf8]{inputenc}
\usepackage[english]{babel}
\usepackage{amsmath}
\usepackage{amsfonts}
\usepackage{amssymb}
\usepackage{mathtools}
\usepackage{amsthm}
\usepackage{float}
\usepackage{etoolbox}
\usepackage{bbm}
\usepackage{hyperref}
\usepackage{graphicx}
\usepackage{xcolor}
\usepackage[T1]{fontenc}
\newtheorem{theorem}{Theorem}
\newtheorem{lemma}{Lemma}
\newtheorem{remark}{Remark}

\newtheorem{corollary}{Corollary}
\newtheorem{assumption}{Assumption}
\usepackage[left=2cm,right=2cm,top=2cm,bottom=2cm]{geometry}

\newcommand{\probability}{\operatorname{\mathbb{P}}\probarg}
\DeclarePairedDelimiterX{\probarg}[1]{(}{)}{%
  \ifnum\currentgrouptype=16 \else\begingroup\fi
  \activatebar#1
  \ifnum\currentgrouptype=16 \else\endgroup\fi
}
\newcommand{\activatebar}{%
  \begingroup\lccode`\~=`\|
  \lowercase{\endgroup\let~}\innermid 
  \mathcode`|=\string"8000
}
\newcommand{\LimitN}{\overset{N\to\infty}{\longrightarrow}}
\newcommand{\LimitD}{\overset{d}{\longrightarrow}}
\newcommand{\LimitP}{\overset{\mathbb{P}}{\longrightarrow}}

\DeclareRobustCommand{\VAN}[3]{#2}
\usepackage[utf8]{inputenc}

\title{Extreme values for the waiting time in large fork-join queues}
\author{Dennis Schol, Maria Vlasiou, and Bert Zwart}
\date{\today}

\begin{document}

\maketitle
\begin{abstract}
   We prove that the scaled maximum steady-state waiting time and the scaled maximum steady-state queue length among $N$ $GI/GI/1$-queues in the $N$-server fork-join queue, converge to a normally distributed random variable as $N\to\infty$. The maximum steady-state waiting time in this queueing system scales around $\frac{1}{\gamma}\log N$, where $\gamma$ is determined by the cumulant generating function $\Lambda$ of the service distribution and solves the Cram\'er-Lundberg equation with stochastic service times and deterministic inter-arrival times. This value $\frac{1}{\gamma}\log N$ is reached at a certain hitting time. The number of arrivals until that hitting time satisfies the central limit theorem, with standard deviation $\frac{\sigma_A}{\sqrt{\Lambda'(\gamma)\gamma}}$. By using distributional Little's law, we can extend this result to the maximum queue length. Finally, we extend these results to a fork-join queue with different classes of servers. 
\end{abstract}
\section{Introduction}
Fork-join queues are a useful modeling tool for congestion in complex networks, such as assembly systems, communication networks, and supply chains. Such networks can be large and assembly is only possible upon availability of all parts. Thus, the bottleneck of the system is caused by the slowest production line in the system. This setting motivates us to investigate such delays in a stylized version of a large fork-join queueing system. In this setting, a key quantity of interest is the behavior of the longest queue when the system is in the steady-state situation. Furthermore, we assume that arrival and service processes are general and mutually independent.

 As we try to model systems with many servers, we are typically interested in the behavior of this random variable as $N\to\infty$. In \cite{MeijerSchol2021}, it is shown that $\max_{i\leq N}(B_i(s)+B_A(s)-\beta s)$ is in the domain of attraction of the normal distribution:
\begin{align}\label{eq: central limit convergence}
    \probability*{\max_{i\leq N}(B_i(s)+B_A(s)-\beta s)>\frac{\sigma^2}{2\beta}\log N+x\sqrt{\log N}}\LimitN\probability*{\frac{\sigma\sigma_A}{\sqrt{2}\beta}X>x},
\end{align}
with $X\overset{d}{=}\mathcal{N}(0,1)$, where $\{B_i(t),t\geq 0\}$ and $\{B_A(t),t\geq 0\}$ are Brownian motions with standard deviations $\sigma$ and $\sigma_A$, respectively. We see from the limit in \eqref{eq: central limit convergence} that $\max_{i\leq N}(B_i(s)+B_A(s)-\beta s)$ centers around $\frac{\sigma^2}{2\beta}\log N$ and deviates with order $\sqrt{\log N}$.

This convergence result provides a prediction of the typical delay. In this study, we aim to extend this result to a more general setting. In particular, we investigate the maximum steady-state waiting time among the $N$ servers with a common arrival process $\max_{i\leq N}W_i(\infty)=\max_{i\leq N}\sup_{k\geq 0}\sum_{j=1}^k(S_i(j)-A(j))$. This expression follows from Lindley's recursion. Furthermore, we have that both $(S_i(j),j\geq 1,1\leq i\leq N)$ and $(A(j),j\geq 1)$ are i.i.d.\ and the inter-arrival times and service times are mutually independent. Thus $S_i(j)$ indicates the service time of the $j$-th customer in queue $i$, $A(j)$ indicates the inter-arrival time between the $(j-1)$-st and the $j$-th customer. We see that the maximum steady-state wating time is a maximum of $N$ dependent random variables, due to the common arrival process $(A(j),j\geq 1)$.

The earliest literature on fork-join queues focuses on systems with two service stations. Analytic results, such as asymptotics on limiting distributions, can be found in \cite{baccelli1985two, flatto1984two, de1988fredholm, wright1992two}. However, due to the complexity of fork-join queues, these results cannot be expanded to fork-join queues with more than two service stations. Thus, most of the work on fork-join queues with more than two service stations is focused on finding approximations of performance measures. For example, an approximation of the distribution of the response time in $M/M/s$ fork-join queues is given in Ko and Serfozo \cite{ko2004response}. Upper and lower bounds for the mean response time of servers, and other performance measures, are given by Nelson, Tantawi \cite{nelson1988approximate} and Baccelli, Makowski \cite{baccelli1989queueing}. These bounds can be used for fork-join queues with large size, but apart from this, there is not much literature on the convergence of the longest queue length in a fork-join queue as $N\to\infty$. Some results can be found in \cite{MeijerSchol2021,scholMOR,schol2022tail}. In \cite{MeijerSchol2021}, the same convergence results are given as in this paper, but only for the Brownian fork-join queue, thus this paper extends on these results.

This paper is organized as follows. In Section \ref{sec: model}, we present our main results; in Theorem \ref{thm: convergence maximum waiting time} we state that the longest steady-state waiting time satisfies a central limit result; in Theorem \ref{thm: convergence maximum queue length} we show that a similar result holds for the longest queue length, and in Corollary \ref{cor: K classes} we present a similar result when the service distributions can differ among the different queues. In Section \ref{sec: heuristic analysis} we give an intuition why the results hold and how we prove these. Section \ref{sec: proofs} is devoted to proofs.

\section{Model}\label{sec: model}
We investigate a fork-join queue with $N$ servers. Each of the $N$ servers has the same arrival stream of jobs and works independently from all other servers but with the same service distribution. In this section, we state the main result for the longest steady-state waiting time in Theorem \ref{thm: convergence maximum waiting time}. We also show that a similar result holds for the maximum queue length in Lemma \ref{lem: dist little law} and Theorem \ref{thm: convergence maximum queue length}. Furthermore, we extend the result in Theorem \ref{thm: convergence maximum queue length} to a heterogeneous model in Corollary \ref{cor: K classes}.

We now specify some properties of the service times and interarrival times in this fork-join queueing system. First, the sequence of non-negative random variables $(S_i(j),i\geq 1,j\geq 1)$ are i.i.d.\ with $S_i(j)\sim S$, and $S_i(j)$ indicating the service time of the $j$-th subtask in queue $i$. Furthermore, the sequence of non-negative random variables $(A(j),j\geq 1)$ are i.i.d.\ with $A(j)\sim A$, $\mathbb{E}[A(j)]=1/\lambda$, $\text{Var}(A(j))=\sigma_A^2$, and $A(j)$ indicating the interarrival time between the $(j-1)$-st and the $j$-th task. Finally, we have that $\mathbb{E}[S_i(j)-A(j)]=-\mu$, with $\mu>0$, and $(A(j),j\geq 1)$ and $(S_i(j),i\geq 1,j\geq 1)$ are mutually independent.

We can now write the cumulative distribution function of the longest steady-state waiting time as the cumulative distribution function of the maximum of $N$ all-time suprema of random walks involving the interarrival and service times.
\begin{lemma}\label{lem: waiting in d sup rw}
For the model given in Section \ref{sec: model} with $W_i(1)=0$ for all $i\leq N$, we have that the longest waiting time in steady state satisfies
\begin{align}
\max_{i\leq N}W_i(\infty)\overset{d}{=}\max_{i\leq N}\sup_{k\geq 0}\sum_{j=1}^k(S_i(j)-A(j)).
\end{align}
\end{lemma}
\begin{proof}
Using Lindley's recursion \cite{lindley1952theory}, we can write the waiting time of tasks in front of server $i$ as 

$$
W_i(n)=\sup_{0\leq k\leq n}\sum_{j=k+1}^{n}(S_i(j)-A(j)).
$$
Thus, the longest steady steady-state waiting time satisfies.
$$
\max_{i\leq N}W_i(n)=\max_{i\leq N}\sup_{0\leq k\leq n}\sum_{j=k+1}^{n}(S_i(j)-A(j)).
$$
We have that 
$$
\probability{\max_{i\leq N}W_i(\infty)\geq x}=\lim_{n\to\infty}\probability{\max_{i\leq N}W_i(n)\geq x}.
$$
Because, 
\begin{align}
\max_{i\leq N}W_i(n)\overset{d}{=}\max_{i\leq N}\sup_{0\leq k\leq n}\sum_{j=1}^{k}(S_i(j)-A(j)),
\end{align}
we obtain the lemma by using the monotone convergence theorem.

\end{proof}

In order to be able to prove convergence of the longest steady-state waiting time, we need some additional structure for the service-time distribution. We define 
\begin{align}
\Lambda(\theta):=\log(\mathbb{E}[\exp(\theta(S-1/\lambda)]).
\end{align}
Moreover, we write $\mathcal{D}(\Lambda):=\{\theta: \Lambda(\theta)<\infty\}$ and $\mathcal{D}^{\circ}(\Lambda)$ as the interior of $\mathcal{D}(\Lambda)$.
\begin{assumption}\label{assump: 1}
We assume there exists a $\gamma>0$ such that 
\begin{enumerate}
    \item $\Lambda(\gamma)=0$,
    \item $\gamma\in \mathcal{D}^{\circ}(\Lambda)$.
\end{enumerate}
\end{assumption}
The first assumption indicates that the random variable $S-1/\lambda$ has a tail that is bounded by an exponential. The second assumption is needed for our proofs. In \cite[Ex. 2.2.24]{dembo2009large}, it is namely stated that when $\gamma\in \mathcal{D}^{\circ}(\Lambda)$, $\Lambda$ is infinitely differentiable at the point $\gamma$. For example, when $S-1/\lambda$ has density function $f_{S-1/\lambda}(x)=c_1\exp(-x)/(1+x^2)$ for $x>0$, where $c_1,\lambda$ are chosen such that $\probability{S-1/\lambda<x}$ is a cumulative distribution function and $\gamma=1$, then the first assumption is satisfied but the second is not, since $\Lambda(\theta)$ is not differentiable at $\theta=\gamma$.
Our main result is given in Theorem \ref{thm: convergence maximum waiting time}.
\begin{theorem}\label{thm: convergence maximum waiting time}
For the model in Section \ref{sec: model} where the sequence of service times $(S_i(j),i\geq 1,j\geq 1)$ satisfies Assumption \ref{assump: 1}, we have that
\begin{align}\label{eq: limit result maximum waiting time}
\frac{\max_{i\leq N}W_i(\infty)-\frac{1}{\gamma}\log N}{\sqrt{\log N}}\LimitD \frac{\sigma_A}{\sqrt{\Lambda'(\gamma)\gamma}}X,
\end{align}
with $X\sim\mathcal{N}(0,1)$, as $N\to\infty$.
\end{theorem}
\begin{lemma}[Distributional Little's Law]\label{lem: dist little law}
Let for $t\geq 0$, $\mathbf{N}_A(t)$ indicate the number of arrivals up to time $t$, where the interarrival times are i.i.d.\ with $A(j)\sim A$. Then
\begin{align}
    \max_{i\leq N}Q_i(\infty)\overset{d}{=}\mathbf{N}_A\left(\max_{i\leq N}W_i(\infty)\right).
\end{align}
\end{lemma}
\begin{proof}
In \cite{haji1971relation}, a short proof is given that for the $GI/GI/1$ queue under the FCFS policy, $Q\overset{d}{=}\mathbf{N}_A(W)$.
We follow the same steps to prove that $\max_{i\leq N}Q_i(\infty)\overset{d}{=}\mathbf{N}_A\left(\max_{i\leq N}W_i(\infty)\right)$.
\par
First, let $t>0$ be given such that the system is in steady state. Furthermore, let $\tilde{W}_i(j)$ be the waiting time of the $i$-th subtask of the $j$-th task numbered backward in time, beginning at time $t$. Thus, $\tilde{W}_i(1)$ is the waiting time of the $i$-th subtask of the last task arriving before time $t$. Now, let the random variable $T(j)$ be such that $t-T(j)$ is the arrival time of the $j$-th task numbered backward in time. Then, observe that the event $\{\max_{i\leq N}Q_i(t)\geq j\}$ is equivalent to the event that at least one subtask of the $j$-th task numbered backward in time is still in the queue at time $t$. Thus, 
$$
\left\{\max_{i\leq N}Q_i(t)\geq j\right\}=\left\{\max_{i\leq N}\tilde{W}_i(j)\geq T(j)\right\},
$$
for $j\geq 1$. The event $\{T(j)\leq x\}$ is equivalent to the event that the number of arrivals during the period $[t-x,t)$ is larger than or equal to $j$. The arrival process is a stationary process, thus the event $\{T(j)\leq x\}$ is equivalent to the event $\{\mathbf{N}_A(x)\geq j\}$. Additionally, the random variables $\max_{i\leq N}\tilde{W}_i(j)$ and $T(j)$ are independent. Therefore, 
$$
\left\{\max_{i\leq N}Q_i(t)\geq j\right\}=\left\{\mathbf{N}_A\big(\max_{i\leq N}\tilde{W}_i(j)\big)\geq j\right\}.
$$
As the system is in steady state, we get that
$$
\left\{\max_{i\leq N}Q_i(\infty)\geq j\right\}=\left\{\mathbf{N}_A\big(\max_{i\leq N}\tilde{W}_i(\infty)\big)\geq j\right\}.
$$
\end{proof}
Now, combining the result in Lemma \ref{lem: dist little law} with the main result in Theorem \ref{thm: convergence maximum waiting time}, we can find a similar convergence result for the maximum queue length in steady state.
\begin{theorem}\label{thm: convergence maximum queue length}
For the model in Section \ref{sec: model} where the sequence of service times $(S_i(j),i\geq 1,j\geq 1)$ satisfies Assumption \ref{assump: 1}, we have that
\begin{align}
\frac{\max_{i\leq N}Q_i(\infty)-\frac{\lambda}{\gamma}\log N}{\sqrt{\log N}}\LimitD \sqrt{\frac{\lambda^2\sigma_A^2}{\Lambda'(\gamma)\gamma}+\frac{\lambda^3\sigma_A^2}{\gamma}}X,
\end{align}
with $X\sim\mathcal{N}(0,1)$, as $N\to\infty$.
\end{theorem}
\begin{proof}
Let $\hat{A}(j)\sim A$, let $(\hat{A}(j),j\geq 1)$ be mutually independent, and $\hat{A}(j)$ and $\max_{i\leq N}W_i(\infty)$ be mutually independent for all $j\geq 1$. Then, using Lemma \ref{lem: dist little law} and Theorem \ref{thm: convergence maximum waiting time}, we get that
\begin{align*}
  &\probability*{\max_{i\leq N}Q_i(\infty)\leq \frac{\lambda}{\gamma}\log N+x\sqrt{\log N}}\\
  &\quad=\probability*{\mathbf{N}_A\left(\max_{i\leq N}W_i(\infty)\right)\leq \big\lfloor \frac{\lambda}{\gamma}\log N+x\sqrt{\log N}\big\rfloor}\\
   &\quad=\probability*{\max_{i\leq N}W_i(\infty)\leq \sum_{j=1}^{\big\lfloor \frac{\lambda}{\gamma}\log N+x\sqrt{\log N}\big\rfloor}\hat{A}(j)}\\
    &\quad=\probability*{\frac{\max_{i\leq N}W_i(\infty)-\frac{1}{\gamma}\log N}{\sqrt{\log N}}\leq \frac{\sum_{j=1}^{\big\lfloor \frac{\lambda}{\gamma}\log N+x\sqrt{\log N}\big\rfloor}\hat{A}(j)-\frac{1}{\gamma}\log N}{\sqrt{\log N}}}\\
    &\quad\LimitN\probability*{\frac{\sigma_A}{\sqrt{\Lambda'(\gamma)\gamma}}X_1\leq \frac{\sigma_A\sqrt{\lambda}}{\sqrt{\gamma}}X_2+ \frac{x}{\lambda}},
\end{align*}
with $X_1,X_2$ independent and standard normally distributed, this convergence holds, as $(\hat{A}(j),j\geq 1)$ and $\max_{i\leq N}W_i(\infty)$ are independent. Thus, the theorem follows.
\end{proof}

Until now, we considered the fork-join queueing system where each server has the same service distribution. In Corollary \ref{cor: K classes}, we show that we can extend the convergence of the longest steady-state waiting time to a more heterogeneous setting. We examine a fork-join queueing system with $N$ servers, where each of these $N$ servers belongs to one of $K$ classes. Additionally, we assume that the size of class $k$ with $k\in\{1,\ldots,K\}$ grows as $\alpha_kN$, as $N$ becomes large, with $0<\alpha_k<1$.  
\begin{corollary}\label{cor: K classes}
Let $K\in\mathbb{N}$, let $k=1,\ldots,K$, furthermore, take an increasing sequence of integers given by $M_0^{(N)},M_1^{(N)},M_2^{(N)},\ldots, M_K^{(N)}>0$ with $M_0^{(N)}=1$, $M_K^{(N)}=N$, and $M_{k}^{(N)}-M_{k-1}^{(N)}\in \mathbb{N}$. Moreover, $(M_{k}^{(N)}-M_{k-1}^{(N)})/N\LimitN \alpha_k\in(0,1]$ with $\sum_{k=1}^K\alpha_k=1$. Let $(S_i(j),j\geq 1,M_{k-1}^{(N)}< i\leq M_{k}^{(N)})$ be i.i.d.\ with $S_i(j)\sim S_k$, $(A(j),j\geq 1)$ be i.i.d.\ with $A(j)\sim A$, $\mathbb{E}[A(j)]=1/\lambda$, $\text{Var}(A(j))=\sigma_A^2$, $\mathbb{E}[S_i(j)-A(j)]=-\mu_k$ with $\mu_k>0$, $\Lambda_k(\theta)=\log(\mathbb{E}[\exp(\theta(S_k-1/\lambda)])$, $\Lambda_k$ satisfies Assumption \ref{assump: 1}. Furthermore, $S_{i_1}(j_1)$ and $S_{i_2}(j_2)$ are mutually independent for all $i_1,i_2,j_1,j_2$. Let $K^*=\arg \min\{\gamma_k,k=1,\ldots,K\}$. We assume that $|K^*|=1$ and $k^*\in K^*$. Then,
\begin{align}\label{eq: K classes limit}
\frac{\max_{i\leq N}W_i(\infty)-\frac{1}{\gamma_{k^*}}\log N}{\sqrt{\log N}}\LimitD \frac{\sigma_A}{\sqrt{\Lambda_{k^*}'(\gamma_{k^*})\gamma_{k^*}}}X,
\end{align}
with $X\sim\mathcal{N}(0,1)$, as $N\to\infty$.
\end{corollary}
\begin{proof}
We prove this corollary by giving an asymptotically sharp lower and upper bound. First, observe that
$$
\max_{i\leq N}W_i(\infty)\geq_{st.} \max_{M_{k^*-1}^{(N)}< i\leq M_{k^*}^{(N)}}\sup_{k\geq 0}\sum_{j=1}^k(S_i(j)-A(j)),
$$
with $X\geq_{st.}Y$ meaning that $\probability{X\geq x}\geq \probability{Y\geq x}$ for all $x$.
Applying the result from Theorem \ref{thm: convergence maximum waiting time} on the lower bound results in \eqref{eq: K classes limit}. By using the union bound we get the following upper bound:
\begin{multline*}
    \probability*{\max_{i\leq N}W_i(\infty)\geq\frac{1}{\gamma_{k^*}}\log N+x\sqrt{\log N}}\\
    =\sum_{l=1}^K\probability*{\max_{M_{l-1}^{(N)} < i\leq M_l^{(N)}}\sup_{k\geq 0}\sum_{j=1}^k(S_i(j)-A(j))\geq\frac{1}{\gamma_{k^*}}\log N+x\sqrt{\log N}}.
\end{multline*}
When $l\neq k^*$, we get after applying the results from Theorem \ref{thm: convergence maximum waiting time} that 
$$
\probability*{\max_{M_{l-1}^{(N)} < i\leq M_l^{(N)}}\sup_{k\geq 0}\sum_{j=1}^k(S_i(j)-A(j))\geq\frac{1}{\gamma_{l}}\log N+x\sqrt{\log N}}
\LimitN 1-\Phi\left(\frac{\sqrt{\Lambda_l'(\gamma_l)\gamma_l}}{\sigma_A}x\right),
$$
with $\Phi$ the cumulative distribution function of a standard normal random variable.
Because $\gamma_{k^*}<\gamma_l$ we get that
$$
\probability*{\max_{M_{l-1}^{(N)} < i\leq M_l^{(N)}}\sup_{k\geq 0}\sum_{j=1}^k(S_i(j)-A(j))\geq\frac{1}{\gamma_{K^*}}\log N+x\sqrt{\log N}}\LimitN 0.
$$
The corollary follows.
\end{proof}
\begin{remark}
In Corollary \ref{cor: K classes} we assume that $|K^*|=1$. The situation that $|K^*|>1$ follows analogously. Assume for instance that $|K^*|=2$, then we can introduce a new random variable $\tilde{S}$ such that $\tilde{S}_i(j)\sim S_1$ with probability $\alpha$ and $\tilde{S}_i(j)\sim S_2$ with probability $1-\alpha$, such that $\gamma_1=\gamma_2=\gamma_{K^*}$. As $N$ is large enough this fork-join queueing system behaves analogous to the original fork-join queue, and for this system $|K^*|=1$.
\end{remark}

We give the proofs of the convergence of the longest steady-state waiting time in Section \ref{sec: proofs}. First, we give a heuristic explanation of why the convergence result in Theorem \ref{thm: convergence maximum waiting time} is true, and we illustrate the structure of the proof. 
\section{Heuristic analysis}\label{sec: heuristic analysis}
To prove Theorem \ref{thm: convergence maximum waiting time}, we analyze lower and upper bounds of the tail probability of the longest steady-state waiting time among the $N$ servers $\probability{\max_{i\leq N}W_i(\infty)>\frac{1}{\gamma}\log N+x\sqrt{\log N}}$ and we show that these lower and upper bounds converge to the same limit as $N\to\infty$. The longest steady-state waiting time has the form $\max_{i\leq N}W_i(\infty)\overset{d}{=}\sup_{k\geq 0}\max_{i\leq N}\sum_{j=1}^k(S_i(j)-A(j))$. Thus the longest steady-state waiting time is the all-time supremum of the maximum of $N$ random walks. For all processes $(X(t),t\geq 0)$, we have for all $t>0$ 
\begin{align}\label{eq: lower bound}
    \probability*{\sup_{s>0}X(s)>x}\geq \probability*{X(t)>x}.
\end{align}
Furthermore, due to the union bound, we have for all $0<t_1<t_2$ that
\begin{align}\label{eq: upper bound}
    \probability*{\sup_{s>0}X(s)>x} \leq \probability*{\sup_{0<s<t_1}X(s)>x}+\probability*{\sup_{t_1\leq s<t_2}X(s)>x}+\probability*{\sup_{s\geq t_2}X(s)>x}.
\end{align} 
We use these types of lower and upper bounds to prove Theorem \ref{thm: convergence maximum waiting time}. Obviously, not all choices of $t,t_1$, and $t_2$ give sharp bounds. We can however make an educated guess about which choices will give the sharpest bounds. Let us first replace the sequence of random variables $(A(j),j\geq 1)$ with their expectation $1/\lambda$. Thus, we look at a simplified fork-join queue with deterministic arrivals. Because the arrivals are deterministic, the waiting times are mutually independent, and we are able to use standard extreme-value theory. We know from the Cram\'er-Lundberg approximation \cite[Ch.\ XIII, Thm.\ 5.2]{asmussen2003applied} that $\probability{\sup_{k\geq 0}\sum_{j=1}^k(S_i(j)-1/\lambda)>x}\sim C\exp(-\gamma x)$, as $x\to\infty$, with $0<C<1$. Thus, $\probability{\sup_{k\geq 0}\sum_{j=1}^k(S_i(j)-1/\lambda)>\frac{1}{\gamma}\log N}\sim C/N$, as $N\to\infty$. Now we can conclude by using basic extreme-value results; see \cite[Thm.\ 5.4.1, p.\ 188]{de2007extreme}, that
$$
\frac{\max_{i\leq N}\sup_{k\geq 0}\sum_{j=1}^k\left(S_i(j)-\frac{1}{\lambda}\right)}{\log N}\LimitP \frac{1}{\gamma},
$$
as $N\to\infty$. Thus, we know that $\max_{i\leq N}\sup_{k\geq 0}\sum_{j=1}^k(S_i(j)-1/\lambda)$ centers around $\frac{1}{\gamma}\log N$. In order to find suitable lower and upper bounds of the form as given in \eqref{eq: lower bound} and \eqref{eq: upper bound}, we need to estimate the hitting time
$$
\tau^{(N)}:=\inf\left\{k\geq 0:\max_{i\leq N}\sum_{j=1}^k\left(S_i(j)-\frac{1}{\lambda}\right)\geq \frac{1}{\gamma}\log N\right\}.
$$
As mentioned before, we have that $\probability{\sup_{k\geq 0}\sum_{j=1}^k\left(S_i(j)-\frac{1}{\lambda}\right)>\frac{1}{\gamma}\log N}\sim C/N$ as $N\to\infty$. Thus, a good estimate $\hat{\tau}^{(N)}$ for $\tau^{(N)}$ should also satisfy the property that
\begin{align}\label{eq: liminf tail prob}
\liminf_{N\to\infty}N\probability*{\sum_{j=1}^{\hat{\tau}^{(N)}}\left(S_i(j)-\frac{1}{\lambda}\right)>\frac{1}{\gamma}\log N}>0    
\end{align}
and
\begin{align}\label{eq: limsup tail prob}
    \limsup_{N\to\infty}N\probability*{\sum_{j=1}^{\hat{\tau}^{(N)}}\left(S_i(j)-\frac{1}{\lambda}\right)>\frac{1}{\gamma}\log N}<\infty.
\end{align}
Now, by using Cram\'er's theorem and by using the fact that $\Lambda$ is at least twice differentiable at $\gamma$, we know that
\begin{align}\label{eq: cramers theorem}
    \lim_{n\to\infty}\frac{1}{n}\log\left(\probability*{\sum_{j=1}^n\left(S_i(j)-\frac{1}{\lambda}\right)\geq nx}\right)= -\Lambda^*(x),
\end{align}
for all $x>\mathbb{E}[S_i(j)-1/\lambda]$ with $\Lambda^*(x)=\sup_{t\in\mathbb{R}}(tx-\Lambda(t))$; see \cite[Ch. XIII, Thm.\ 2.1 (2.3)]{asmussen2003applied}. We write $\hat{\tau}^{(N)}=\hat{c}\log N$. Then we can conclude from Equation \eqref{eq: cramers theorem} that
\begin{align}\label{eq: cramers theorem2}
    \lim_{N\to\infty}\frac{1}{\log N}\log\left(\probability*{\sum_{j=1}^{\lfloor\hat{c}\log N\rfloor}\left(S_i(j)-\frac{1}{\lambda}\right)\geq x\hat{c}\log N}\right)=-\Lambda^*(x)\hat{c}.
\end{align}
Thus, in order to find a good estimate $\hat{\tau}^{(N)}$ for the hitting time $\tau^{(N)}$ we need to solve two equations. First, $x\hat{c}=1/\gamma$, because we know that the longest steady-state waiting time under deterministic arrivals is approximately equal to $\frac{1}{\gamma}\log N$. Therefore the expression $x\hat{c}\log N$ in \eqref{eq: cramers theorem2} should be the same as $\frac{1}{\gamma}\log N$. Second, $-\Lambda^*(x)\hat{c}=-1$, because we know from \eqref{eq: liminf tail prob}, \eqref{eq: limsup tail prob}, and \eqref{eq: cramers theorem2} that for large $N$
$$
\probability*{\sum_{j=1}^{\lfloor\hat{c}\log N\rfloor}\left(S_i(j)-\frac{1}{\lambda}\right)\geq x\hat{c}\log N}\approx \frac{1}{N}=\exp(-\Lambda^*(x)\hat{c}\log N).
$$
Combining these two equations gives $\hat{c}=\frac{1}{\Lambda'(\gamma)\gamma}$ and $x=\Lambda'(\gamma)$. Clearly, $x\hat{c}=1/\gamma$, and 
$$
\Lambda^*(x)\hat{c}=\frac{\Lambda^*(\Lambda'(\gamma))}{\gamma\Lambda'(\gamma)}.
$$
From \cite[Lem.\ 2.2.5(c)]{dembo2009large}, we know that $\Lambda^*(\Lambda'(\gamma))=\gamma\Lambda'(\gamma)$, thus indeed, $\Lambda^*(x)\hat{c}=1$. Finally, we can conclude that $\hat{\tau}^{(N)}=\hat{c}\log N=\frac{1}{\gamma\Lambda'(\gamma)}\log N$. Obviously, in order to be a good estimation for a hitting time we need to have that $\Lambda'(\gamma)>0$. This is the case because $\Lambda(\theta)$ is convex; see \cite[Ch. XIII, Thm.\ 5.1]{asmussen2003applied}.

Until this point, we know the first-order scaling of the largest of $N$ steady-state waiting times with deterministic arrivals, and we can give an estimation of the hitting time of this value. Now, we can use these results to obtain a second-order convergence result for the longest steady-state waiting time with stochastic arrivals. Following the analysis above together with the lower bound in \eqref{eq: lower bound}, we see that
\begin{multline}\label{eq: lower bound2}
\probability*{\frac{\max_{i\leq N}W_i(\infty)-\frac{1}{\gamma}\log N}{\sqrt{\log N}}\geq x}\\\geq \probability*{\frac{\max_{i\leq N}\sup_{\big(\frac{1}{\Lambda'(\gamma)\gamma}-\epsilon\big)\log N<k<\frac{1}{\Lambda'(\gamma)\gamma}\log N}\sum_{j=1}^{k}(S_i(j)-A(j))-\frac{1}{\gamma}\log N}{\sqrt{\log N}}\geq x},
\end{multline}
with $\epsilon>0$ and small.
\iffalse
For this lower bound, we can write 
\begin{multline}
\frac{\max_{i\leq N}\sup_{\big(\frac{1}{\Lambda'(\gamma)\gamma}-\epsilon\big)\log N<k<\frac{1}{\Lambda'(\gamma)\gamma}\log N}\sum_{j=1}^{k}(S_i(j)-A(j))-\frac{1}{\gamma}\log N}{\sqrt{\log N}}\\
\geq \frac{\max_{i\leq N}\sum_{j=1}^{\big\lfloor\frac{1}{\Lambda'(\gamma)\gamma}\log N\big\rfloor}\left(S_i(j)-\frac{1}{\lambda}\right)-\frac{1}{\gamma}\log N}{\sqrt{\log N}}+\frac{\sum_{j=1}^{\big\lfloor\frac{1}{\Lambda'(\gamma)\gamma}\log N\big\rfloor}(\frac{1}{\lambda}-A(j))}{\sqrt{\log N}}. 
\end{multline}
Obviously, the second term on the right-hand side of \eqref{eq: lower bound2} satisfies the central limit theorem.\fi In Lemma \ref{lem: lower bound}, we prove that the right-hand side in \eqref{eq: lower bound2} converges to a function that is close to the tail probability of a normally distributed random variable. Furthermore, we show in Lemmas \ref{lem: supremum 0 t}, \ref{lem: t-e t+e}, and \ref{lem: t+e infty}, that this lower bound is sharp. To achieve this, we first divide the supremum over all positive numbers in the random variable $\max_{i\leq N}W_i(\infty)$ in three parts. After that, we take the supremum over the intervals $\left[0,\left(\frac{1}{\Lambda'(\gamma)\gamma}-\epsilon\right)\log N\right]$, $\left(\left(\frac{1}{\Lambda'(\gamma)\gamma}-\epsilon\right)\log N,\left(\frac{1}{\Lambda'(\gamma)\gamma}+\epsilon\right)\log N\right]$, and $\left(\left(\frac{1}{\Lambda'(\gamma)\gamma}+\epsilon\right)\log N,\infty\right)$, with $\epsilon>0$ and small. Consequently, we show that the tail probabilities of the first and third suprema of the maximum of $N$ random walks asymptotically vanish, while 
$$
\probability*{\max_{i\leq N}\sup_{\big(\frac{1}{\Lambda'(\gamma)\gamma}-\epsilon\big)\log N<k<\big(\frac{1}{\Lambda'(\gamma)\gamma}+\epsilon\big)\log N}\sum_{j=1}^k(S_i(j)-A(j))>\frac{1}{\gamma}\log N+x\sqrt{\log N}}
$$
converges to a limit close to the lower bound as $N\to\infty$. 
\begin{remark}
The lower bound presented in Equation \eqref{eq: lower bound} gives us information about the convergence rate of the result in Theorem \ref{thm: convergence maximum waiting time}. From the Berry-Ess\'een theorem \cite{michel1981constant}, we know that when $\frac{1}{\sqrt{n}}\sum_{i=1}^nX_i\LimitD X\sim\mathcal{N}(0,1)$, the convergence rate is of order $1/\sqrt{n}$. Thus, the lower bound in \eqref{eq: lower bound} shows that the convergence rate is of order $1/\sqrt{\log N}$.
\end{remark}
\section{Proofs}\label{sec: proofs}
\begin{lemma}\label{lem: lower bound}
Given the model in Section \ref{sec: model} where the sequence of service times $(S_i(j),i\geq 1,j\geq 1)$ satisfies Assumption \ref{assump: 1}, $0<\epsilon<\frac{1}{\Lambda'(\gamma)\gamma}$, $t_1^{(N)}=\left(\frac{1}{\Lambda'(\gamma)\gamma}-\epsilon\right)\log N$, and $t_2^{(N)}=\frac{1}{\Lambda'(\gamma)\gamma}\log N$, then for all $x\in \mathbb{R}$, we have that
\begin{multline}\label{eq: lower bound convergence}
\liminf_{N\to\infty}\probability*{\max_{i\leq N}\sup_{t_1^{(N)}<k<t_2^{(N)}}\sum_{j=1}^k(S_i(j)-A(j))>\frac{1}{\gamma}\log N+x\sqrt{\log N}}\\
\geq  \probability*{\sigma_A\sqrt{\frac{1}{\Lambda'(\gamma)\gamma}-\epsilon}X_1-\sigma_A\sqrt{\epsilon}\left|X_2\right|>x},
\end{multline}
with $X_1,X_2\sim\mathcal{N}(0,1)$ and independent.

\end{lemma}
\begin{proof}
In order to prove this convergence result, we first bound 

$$
\max_{i\leq N}\sup_{t_1^{(N)}<k<t_2^{(N)}}\sum_{j=1}^k(S_i(j)-A(j))
\geq\max_{i\leq N}\sup_{t_1^{(N)}<k<t_2^{(N)}}\sum_{j=1}^k\bigg(S_i(j)-\frac{1}{\lambda}\bigg)
+\inf_{t_1^{(N)}< k< t_2^{(N)}}\sum_{j=1}^{k}\bigg(\frac{1}{\lambda}-A(j)\bigg).
$$
We treat the terms on the right-hand side separately. We first prove that
\begin{align}\label{eq: lower bound arrival process clt}
\frac{\inf_{t_1^{(N)}< k< t_2^{(N)}}\sum_{j=1}^{k}\left(\frac{1}{\lambda}-A(j)\right)}{\sqrt{\log N}}\LimitD \sigma_A\sqrt{\frac{1}{\Lambda'(\gamma)\gamma}-\epsilon}X_1-\sigma_A\sqrt{\epsilon}\left|X_2\right|,
\end{align}
as $N\to\infty$.
Afterwards, we prove that 
\begin{align}\label{eq: independent part to zero}
\frac{\max_{i\leq N}\sup_{t_1^{(N)}<k<t_2^{(N)}}\sum_{j=1}^k\left(S_i(j)-\frac{1}{\lambda}\right)-\frac{1}{\gamma}\log N}{\sqrt{\log N}}\LimitP 0,
\end{align} 
as $N\to\infty$.

The first convergence result follows from Donsker's theorem. The left-hand side in \eqref{eq: lower bound arrival process clt} is an infimum of a random walk with drift 0. Then for $(B(t),t\geq 0)$ a Brownian motion with drift 0 and standard deviation 1, by using Donsker's theorem \cite{donsker1951invariance} and the fact that the infimum is a continuous functional, we obtain that
$$
\probability*{\frac{\inf_{t_1^{(N)}< k< t_2^{(N)}}\sum_{j=1}^{k}(\frac{1}{\lambda}-A(j))}{\sqrt{\log N}}>x}\LimitN \probability*{\inf_{\left(\frac{1}{\Lambda'(\gamma)\gamma}-\epsilon\right)<s<\frac{1}{\Lambda'(\gamma)\gamma}}\sigma_A B(s)>x}.
$$
Furthermore, we can rewrite $$\inf_{\frac{1}{\Lambda'(\gamma)\gamma}-\epsilon<s<\frac{1}{\Lambda'(\gamma)\gamma}}\sigma_A B(s)\overset{d}{=}\sigma_AB\left(\frac{1}{\Lambda'(\gamma)\gamma}-\epsilon\right)-\inf_{0<s<\epsilon}\sigma_A\tilde{B}(s),$$ where $\tilde{B}$ is an independent copy of $B$. Obviously, we have that $\sigma_AB\left(\frac{1}{\Lambda'(\gamma)\gamma}-\epsilon\right)\overset{d}{=}\sigma_A\sqrt{\frac{1}{\Lambda'(\gamma)\gamma}-\epsilon}X_1$ with $X_1\sim\mathcal{N}(0,1)$. Because
$\inf_{0<s<\epsilon}\sigma_A\tilde{B}(s)\overset{d}{=}\sigma_A\sqrt{\epsilon}|X_2|$, with $X_2\sim\mathcal{N}(0,1)$, we have that the limit in \eqref{eq: lower bound arrival process clt} follows.

In order to prove the second convergence result, we define for $A\in\mathcal{F}_k$, with $\{\mathcal{F}_k,k\geq 1\}$ the natural filtration, the probability measure 
$$
\mathbb{P}_i(A):=\mathbb{E}\left[\exp\left(\gamma\sum_{j=1}^k\left(S_i(j)-\frac{1}{\lambda}\right)\right)\mathbbm{1}(A)\right];
$$
see \cite[Ch. XIII, Par. 3]{asmussen2003applied}. Now, we know that 
$$
\mathbb{E}_i\left[S_i(j)-\frac{1}{\lambda}\right]=\mathbb{E}\left[\left(S_i(j)-\frac{1}{\lambda}\right)\exp\left(\gamma\left(S_i(j)-\frac{1}{\lambda}\right)\right)\right]=\Lambda'(\gamma).
$$
Thus, by checking the conditions in \cite[Ch.\ XIII, Thm.\ 5.6]{asmussen2003applied}, we see that 
\begin{multline}\label{eq: supremum conv 1}
\probability*{\sup_{0\leq k< t_2^{(N)}}\sum_{j=1}^{k}\left(S_i(j)-\frac{1}{\lambda}\right)\geq\frac{1}{\gamma}\log N+x\sqrt{\log N}}\\
=C\exp\bigg(-\gamma\bigg(\frac{1}{\gamma}\log N+x\sqrt{\log N}\bigg)\bigg)\Phi\bigg(-x\frac{\sqrt{\gamma\Lambda'(\gamma)}}{\sqrt{\Lambda''(\gamma)}}\bigg)(1+o(1)).
\end{multline}
With the same approach, we get from \cite[Ch.\ XIII, Thm.\ 5.6]{asmussen2003applied} that
\begin{align}\label{eq: supremum conv 2}
\probability*{\sup_{0\leq k< t_1^{(N)}}\sum_{j=1}^{k}\left(S_i(j)-\frac{1}{\lambda}\right)\geq\frac{1}{\gamma}\log N+x\sqrt{\log N}}
=o\bigg(C\exp\bigg(-\gamma\bigg(\frac{1}{\gamma}\log N+x\sqrt{\log N}\bigg)\bigg)\bigg),
\end{align}
as $N\to\infty$, for all $x\in\mathbb{R}$. By applying the union bound, we get that 
\begin{align*}
&\probability*{\sup_{0\leq k< t_2^{(N)}}\sum_{j=1}^{k}\left(S_i(j)-\frac{1}{\lambda}\right)\geq\frac{1}{\gamma}\log N+x\sqrt{\log N}}\\
&\quad\leq\probability*{\sup_{0\leq k< t_1^{(N)}}\sum_{j=1}^{k}\left(S_i(j)-\frac{1}{\lambda}\right)\geq\frac{1}{\gamma}\log N+x\sqrt{\log N}}+\probability*{\sup_{t_1^{(N)}< k< t_2^{(N)}}\sum_{j=1}^{k}\left(S_i(j)-\frac{1}{\lambda}\right)\geq\frac{1}{\gamma}\log N+x\sqrt{\log N}}\\
&\quad\leq\probability*{\sup_{0\leq k< t_1^{(N)}}\sum_{j=1}^{k}\left(S_i(j)-\frac{1}{\lambda}\right)\geq\frac{1}{\gamma}\log N+x\sqrt{\log N}}+\probability*{\sup_{0\leq k< t_2^{(N)}}\sum_{j=1}^{k}\left(S_i(j)-\frac{1}{\lambda}\right)\geq\frac{1}{\gamma}\log N+x\sqrt{\log N}}.
\end{align*}
We can conclude from these bounds, together with \eqref{eq: supremum conv 1} and \eqref{eq: supremum conv 2} that
\begin{multline}
\probability*{\sup_{t_1^{(N)}< k< t_2^{(N)}}\sum_{j=1}^{k}\left(S_i(j)-\frac{1}{\lambda}\right)\geq\frac{1}{\gamma}\log N+x\sqrt{\log N}}\\
=C\exp\bigg(-\gamma\bigg(\frac{1}{\gamma}\log N+x\sqrt{\log N}\bigg)\bigg)\Phi\bigg(-x\frac{\sqrt{\gamma\Lambda'(\gamma)}}{\sqrt{\Lambda''(\gamma)}}\bigg)(1+o(1)).
\end{multline}
By using this expression it is easy to derive that for $x>0$
\begin{multline*}
\probability*{\max_{i\leq N}\sup_{t_1^{(N)}< k<t_2^{(N)}}\sum_{j=1}^{k}\left(S_i(j)-\frac{1}{\lambda}\right)\leq\frac{1}{\gamma}\log N+x\sqrt{\log N}}\\
=\probability*{\sup_{t_1^{(N)}< k< t_2^{(N)}}\sum_{j=1}^{k}\left(S_i(j)-\frac{1}{\lambda}\right)\leq\frac{1}{\gamma}\log N+x\sqrt{\log N}}^N\LimitN 1.
\end{multline*}
Similarly, for $x<0$,
\begin{multline*}
\probability*{\max_{i\leq N}\sup_{t_1^{(N)}< k<t_2^{(N)}}\sum_{j=1}^{k}\left(S_i(j)-\frac{1}{\lambda}\right)\leq\frac{1}{\gamma}\log N+x\sqrt{\log N}}\\
=\probability*{\sup_{t_1^{(N)}< k< t_2^{(N)}}\sum_{j=1}^{k}\left(S_i(j)-\frac{1}{\lambda}\right)\leq\frac{1}{\gamma}\log N+x\sqrt{\log N}}^N\LimitN 0.
\end{multline*}

Combining these two results gives us the limit in \eqref{eq: independent part to zero}. Finally, the convergence result in \eqref{eq: lower bound convergence} follows from the two limits in \eqref{eq: lower bound arrival process clt} and \eqref{eq: independent part to zero}. 
\end{proof}
\begin{lemma}\label{lem: supremum 0 t}
Given the model in Section \ref{sec: model} where the sequence of service times $(S_i(j),i\geq 1,j\geq 1)$ satisfies Assumption \ref{assump: 1}, $t_1^{(N)}=\left(\frac{1}{\Lambda'(\gamma)\gamma}-\epsilon\right)\log N$, $\delta=\frac{\delta_1}{\Lambda'(\gamma)\gamma}+\delta_2$ with $\delta_{1,2}>0$ and small, and $\epsilon=\delta^{1/4}$, then for all $x\in \mathbb{R}$, we have that
\begin{align}\label{eq: sup 0 t}
\probability*{\max_{i\leq N}\sup_{0\leq k<t_1^{(N)}}\sum_{j=1}^k(S_i(j)-A(j))>\frac{1}{\gamma}\log N+x\sqrt{\log N}}\LimitN 0.
\end{align}
\end{lemma}
\begin{proof}
We derive upper bounds for the left-hand side of \eqref{eq: sup 0 t} that converge to 0 as $N\to\infty$.  

We get by using the subadditivity property of the sup operator and the union bound that
\begin{align}
&\probability*{\max_{i\leq N}\sup_{0\leq k<t_1^{(N)}}\sum_{j=1}^k(S_i(j)-A(j))>\frac{1}{\gamma}\log N}\\
&\quad\leq \probability*{\max_{i\leq N}\sup_{0\leq k<t_1^{(N)}}\sum_{j=1}^k\left(S_i(j)-\frac{1}{\lambda}+\delta_1\right)>\left(\frac{1}{\gamma}-\delta_2\right)\log N}\label{subeq: upper bound 0}\\
&\quad\quad+\probability*{\sup_{k\geq 0}\sum_{j=1}^k\left(\frac{1}{\lambda}-\delta_1-A(j)\right)>\delta_2\log N+x\sqrt{\log N}}\label{subeq: upper bound 1 performance}.
\end{align}
First, because $\mathbb{E}[\frac{1}{\lambda}-\delta_1-A(j)]<0$, we get that 
$$\probability*{\sup_{k\geq 0}\sum_{j=1}^k\left(\frac{1}{\lambda}-\delta_1-A(j)\right)>\delta_2\log N+x\sqrt{\log N}}\LimitN 0.$$
Second, we can bound the term in \eqref{subeq: upper bound 0} as follows;
\begin{multline*}
\probability*{\max_{i\leq N}\sup_{0\leq k<t_1^{(N)}}\sum_{j=1}^k\left(S_i(j)-\frac{1}{\lambda}+\delta_1\right)>\left(\frac{1}{\gamma}-\delta_2\right)\log N}\\
\leq \probability*{\max_{i\leq N}\sup_{0\leq k<t_1^{(N)}}\sum_{j=1}^k\left(S_i(j)-\frac{1}{\lambda}\right)>\left(\frac{1}{\gamma}-\frac{\delta_1}{\Lambda'(\gamma)\gamma}-\delta_2\right)\log N}.
\end{multline*}
Now, we can bound this further to
\begin{multline*}
    \probability*{\max_{i\leq N}\sup_{0\leq k<t_1^{(N)}}\sum_{j=1}^k\left(S_i(j)-\frac{1}{\lambda}\right)>\left(\frac{1}{\gamma}-\frac{\delta_1}{\Lambda'(\gamma)\gamma}-\delta_2\right)\log N}\\
    \leq \sum_{k=0}^{\lfloor t_1^{(N)}\rfloor }N\probability*{\sum_{j=1}^k\left(S_i(j)-\frac{1}{\lambda}\right)>\left(\frac{1}{\gamma}-\frac{\delta_1}{\Lambda'(\gamma)\gamma}-\delta_2\right)\log N}.
\end{multline*}
By using Chernoff's bound we obtain that for $\Lambda(\theta)<\infty$
\begin{align}
 &\sum_{k=0}^{\lfloor t_1^{(N)}\rfloor }N\probability*{\sum_{j=1}^k\left(S_i(j)-\frac{1}{\lambda}\right)>\left(\frac{1}{\gamma}-\frac{\delta_1}{\Lambda'(\gamma)\gamma}-\delta_2\right)\log N}\\
&\quad\leq N \sum_{k=0}^{\lfloor t_1^{(N)}\rfloor }\exp(k \Lambda(\theta))\exp\bigg(-\theta\bigg(\frac{1}{\gamma}-\frac{\delta_1}{\Lambda'(\gamma)\gamma}-\delta_2\bigg)\log N\bigg)\nonumber\\
&\quad =  N\frac{-1+\exp\left((\lfloor t_1^{(N)}\rfloor+1)\Lambda(\theta)\right)}{\exp(\Lambda(\theta))-1}\exp\bigg(-\theta\bigg(\frac{1}{\gamma}-\frac{\delta_1}{\Lambda'(\gamma)\gamma}-\delta_2\bigg)\log N\bigg).\label{subeq: upper bound 2 performance}
\end{align}
Now, 
\begin{multline*}
\frac{\log\left(N\frac{-1+\exp\left((\lfloor t_1^{(N)}\rfloor+1)\Lambda(\theta)\right)}{\exp(\Lambda(\theta))-1}\exp\bigg(-\theta\bigg(\frac{1}{\gamma}-\frac{\delta_1}{\Lambda'(\gamma)\gamma}-\delta_2\bigg)\log N\bigg)\right)}{\log N}\\
\LimitN 1-\left(\theta\left(\frac{1}{\gamma}-\frac{\delta_1}{\Lambda'(\gamma)\gamma}-\delta_2\right)-\left(\frac{1}{\Lambda'(\gamma)\gamma}-\epsilon\right)\Lambda(\theta)\right).
\end{multline*}
In order to make the bound in \eqref{subeq: upper bound 2 performance} as sharp as possible, we need to choose a convenient $\theta$. The choice of $\theta$ that gives the sharpest bound maximizes the function $\theta\left(\frac{1}{\gamma}-\frac{\delta_1}{\Lambda'(\gamma)\gamma}-\delta_2\right)-\left(\frac{1}{\Lambda'(\gamma)\gamma}-\epsilon\right)\Lambda(\theta)$. We have that $\delta=\frac{\delta_1}{\Lambda'(\gamma)\gamma}+\delta_2$ and $\epsilon=\delta^{1/4}$. Furthermore, we choose $\theta=\gamma+\sqrt{\delta}$. This gives us a sharp enough bound in \eqref{subeq: upper bound 2 performance}. We obviously have that
$$
\sup_{\eta\in\mathbb{R}}\left(\eta\left(\frac{1}{\gamma}-\delta\right)-\left(\frac{1}{\Lambda'(\gamma)\gamma}-\delta^{1/4}\right)\Lambda(\eta)\right)
\geq \left((\gamma+\sqrt{\delta})\left(\frac{1}{\gamma}-\delta\right)-\left(\frac{1}{\Lambda'(\gamma)\gamma}-\delta^{1/4}\right)\Lambda(\gamma+\sqrt{\delta})\right).
$$
The first order Taylor series of $\Lambda(\gamma+\sqrt{\delta})$ around $\gamma$ gives
$$
\Lambda(\gamma+\sqrt{\delta})=\Lambda(\gamma)+\sqrt{\delta}\Lambda'(\gamma)+O(\delta)=\sqrt{\delta}\Lambda'(\gamma)+O(\delta).
$$
Thus,
$$
\left((\gamma+\sqrt{\delta})\left(\frac{1}{\gamma}-\delta\right)-\left(\frac{1}{\Lambda'(\gamma)\gamma}-\delta^{1/4}\right)\Lambda(\gamma+\sqrt{\delta})\right)=1+\delta^{3/4}\Lambda'(\gamma)+O(\delta)>1,
$$
for $\delta$ small enough.
Thus the expression in \eqref{subeq: upper bound 2 performance} is upper bounded by the term $N^{-\delta^{3/4}\Lambda'(\gamma)-O(\delta)}\LimitN 0$.
\end{proof}
\begin{lemma}\label{lem: t-e t+e}
Given the model in Section \ref{sec: model} where the sequence of service times $(S_i(j),i\geq 1,j\geq 1)$ satisfies Assumption \ref{assump: 1}, $0<\epsilon<\frac{1}{\Lambda'(\gamma)\gamma}$, $t_1^{(N)}=\left(\frac{1}{\Lambda'(\gamma)\gamma}-\epsilon\right)\log N$, and $t_3^{(N)}=\left(\frac{1}{\Lambda'(\gamma)\gamma}+\epsilon\right)\log N$, then for all $x\in \mathbb{R}$, we have that
\begin{multline}
\limsup_{N\to\infty}\probability*{\max_{i\leq N}\sup_{t_1^{(N)}\leq k<t_3^{(N)}}\sum_{j=1}^k(S_i(j)-A(j))>\frac{1}{\gamma}\log N+x\sqrt{\log N}}\\
\leq  \probability*{\sigma_A\sqrt{\frac{1}{\Lambda'(\gamma)\gamma}-\epsilon}X_1+\sigma_A\sqrt{2\epsilon}\left|X_2\right|>x},
\end{multline}
with $X_1,X_2\sim\mathcal{N}(0,1)$ and independent.
\end{lemma}
\begin{proof}
In order to prove this lemma, we first rewrite 
\begin{align}
   &\frac{\max_{i\leq N}\sup_{t_1^{(N)}\leq k<t_3^{(N)}}\sum_{j=1}^k(S_i(j)-A(j))-\frac{1}{\gamma}\log N}{\sqrt{\log N}}\nonumber\\
    &\quad\leq \frac{\max_{i\leq N}\sup_{t_1^{(N)}\leq k<t_3^{(N)}}\sum_{j=1}^k\left(S_i(j)-\frac{1}{\lambda}\right)-\frac{1}{\gamma}\log N}{\sqrt{\log N}}+ \frac{\sup_{t_1^{(N)}\leq k<t_3^{(N)}}\sum_{j=1}^k\left(\frac{1}{\lambda}-A(j)\right)}{\sqrt{\log N}}
   \nonumber\\
    &\quad\leq \frac{\max_{i\leq N}\sup_{k\geq 0}\sum_{j=1}^k\left(S_i(j)-\frac{1}{\lambda}\right)-\frac{1}{\gamma}\log N}{\sqrt{\log N}}+ \frac{\sup_{t_1^{(N)}\leq k<t_3^{(N)}}\sum_{j=1}^k\left(\frac{1}{\lambda}-A(j)\right)}{\sqrt{\log N}}.\label{subeq: upper bound term 1}
\end{align}
We first look at the first term in \eqref{subeq: upper bound term 1}. This term gives the rescaled longest steady-state waiting time of $N$ i.i.d.\ $D/G/1$ queues. We know that 
$$\probability*{\sup_{k\geq 0}\sum_{j=1}^k\left(S_i(j)-\frac{1}{\lambda}\right)>x}\sim C\exp(-\gamma x),$$
as $x\to\infty$, with $0<C<1$; see \cite[Ch. XIII, Thm.\ 5.2]{asmussen2003applied}. Thus for $x>0$,
\begin{multline*}
\probability*{\frac{\max_{i\leq N}\sup_{k\geq 0}\sum_{j=1}^k\left(S_i(j)-\frac{1}{\lambda}\right)-\frac{1}{\gamma}\log N}{\sqrt{\log N}}>x}\\\sim 1-\left(1-C\exp(-\gamma (1/\gamma\log N+x\sqrt{\log N}))\right)^N\LimitN 0.
\end{multline*}
Similarly, for $x<0$,
\begin{multline*}
\probability*{\frac{\max_{i\leq N}\sup_{k\geq 0}\sum_{j=1}^k\left(S_i(j)-\frac{1}{\lambda}\right)-\frac{1}{\gamma}\log N}{\sqrt{\log N}}>x}\\\sim 1-\left(1-C\exp(-\gamma (1/\gamma\log N+x\sqrt{\log N}))\right)^N\LimitN 1.
\end{multline*}
Thus, the first term in \eqref{subeq: upper bound term 1} converges in probability to 0.

Now, we prove convergence of the tail probability of the second term in \eqref{subeq: upper bound term 1}. This term is a supremum of a random walk with drift 0. Then for $(B(t),t\geq 0)$ a Brownian motion with drift 0 and standard deviation 1, by using Donsker's theorem \cite{donsker1951invariance} and the fact that the supremum is a continuous functional, we obtain with a similar analysis as in Lemma \ref{lem: lower bound}, that
$$
\probability*{\frac{\sup_{t_1^{(N)}\leq k<t_3^{(N)}}\sum_{j=1}^k(\frac{1}{\lambda}-A(j))}{\sqrt{\log N}}>x}
\LimitN \probability*{\sigma_A\sqrt{\frac{1}{\Lambda'(\gamma)\gamma}-\epsilon}X_1+\sigma_A\sqrt{2\epsilon}\left|X_2\right|>x}.
$$
\end{proof}
\begin{lemma}\label{lem: t+e infty}
Given the model in Section \ref{sec: model} where the sequence of service times $(S_i(j),i\geq 1,j\geq 1)$ satisfies Assumption \ref{assump: 1}, $\delta=\frac{\delta_1}{\Lambda'(\gamma)\gamma}+\delta_2$ with $\delta_{1,2}>0$ and small, $\epsilon=\delta^{1/4}$, and $t_3^{(N)}=\left(\frac{1}{\Lambda'(\gamma)\gamma}+\epsilon\right)\log N$, then for all $x\in \mathbb{R}$, we have that
\begin{align}
\probability*{\max_{i\leq N}\sup_{k\geq t_3^{(N)}}\sum_{j=1}^k(S_i(j)-A(j))>\frac{1}{\gamma}\log N+x\sqrt{\log N}}\LimitN 0.
\end{align}
\end{lemma}
\begin{proof}
As in the proof of Lemma \ref{lem: supremum 0 t}, we derive upper bounds for $$\probability*{\max_{i\leq N}\sup_{k\geq t_3^{(N)}}\sum_{j=1}^k(S_i(j)-A(j))>\frac{1}{\gamma}\log N+x\sqrt{\log N}}$$ that converge to 0 as $N\to\infty$.

First, we see that by using subadditivity and the union bound, we obtain
\begin{align*}
&\probability*{\max_{i\leq N}\sup_{k\geq t_3^{(N)}}\sum_{j=1}^k(S_i(j)-A(j))>\frac{1}{\gamma}\log N+x\sqrt{\log N}}\\
&\quad\leq \probability*{\max_{i\leq N}\sup_{k\geq t_3^{(N)}}\sum_{j=1}^k\left(S_i(j)-\frac{1}{\lambda}+\delta_1\right)>\left(\frac{1}{\gamma}-\delta_2\right)\log N}\\
&\quad\quad+\probability*{\sup_{k\geq 0}\sum_{j=1}^k\left(\frac{1}{\lambda}-\delta_1-A(j)\right)>\delta_2\log N+x\sqrt{\log N}}.
\end{align*}
As in the proof of Lemma \ref{lem: supremum 0 t}, we have that 
$$\probability*{\sup_{k\geq 0}\sum_{j=1}^k\left(\frac{1}{\lambda}-\delta_1-A(j)\right)>\delta_2\log N+x\sqrt{\log N}}\LimitN 0.$$
Furthermore, observe that
$\log\mathbb{E}[\exp(\theta(S_i(j)-1/\lambda+\delta_1))]=\Lambda(\theta)+\theta\delta_1$. Now, as in the proof of Lemma \ref{lem: supremum 0 t}, we can bound 
\begin{align}
    &\probability*{\max_{i\leq N}\sup_{k\geq t_3^{(N)}}\sum_{j=1}^k\left(S_i(j)-\frac{1}{\lambda}+\delta_1\right)>\left(\frac{1}{\gamma}-\delta_2\right)\log N}\\
  &\quad   \leq N\sum_{k=\lfloor t_3^{(N)}\rfloor }^{\infty}\probability*{\sum_{j=1}^k\left(S_i(j)-\frac{1}{\lambda}+\delta_1\right)>\left(\frac{1}{\gamma}-\delta_2\right)\log N}\\
&\quad\leq N \sum_{k=\lfloor t_3^{(N)}\rfloor }^{\infty}\exp(k (\Lambda(\theta)+\theta\delta_1))\exp\left(-\theta\left(\frac{1}{\gamma}-\delta_2\right)\log N\right)\nonumber\\
&\quad =  N\frac{\exp\left(\lfloor t_3^{(N)}\rfloor(\Lambda(\theta)+\theta\delta_1)\right)}{\exp(\Lambda(\theta)+\theta\delta_1)-1}\exp\left(-\theta\left(\frac{1}{\gamma}-\delta_2\right)\log N\right),\label{subeq: upper bound 3}
\end{align}
when $\Lambda(\theta)+\theta\delta_1<0$. When $\Lambda(\theta)+\theta\delta_1\geq 0$ the sum in the upper bound diverges to $\infty$.
Now, for the case $\Lambda(\theta)+\theta\delta_1<0$, we have that
\begin{align*}
\frac{\log\left(N\frac{\exp\left(\lfloor t_3^{(N)}\rfloor(\Lambda(\theta)+\theta\delta_1)\right)}{\exp(\Lambda(\theta)+\theta\delta_1)-1}\exp\left(-\theta\left(\frac{1}{\gamma}-\delta_2\right)\log N\right)\right)}{\log N}\LimitN 1+\left(\frac{1}{\Lambda'(\gamma)\gamma}+\epsilon\right)(\Lambda(\theta)+\theta\delta_1)-\theta\left(\frac{1}{\gamma}-\delta_2\right).
\end{align*}
As in the proof of Lemma \ref{lem: supremum 0 t}, we have $\delta=\frac{\delta_1}{\Lambda'(\gamma)\gamma}+\delta_2$ and $\epsilon=\delta^{1/4}$. We now get after a similar derivation as in the proof of Lemma \ref{lem: supremum 0 t} that $\theta=\gamma-\sqrt{\delta}$ gives a sharp bound. First, observe that $\Lambda(\gamma-\sqrt{\delta})=-\sqrt{\delta}\Lambda'(\gamma)+O(\delta)$, thus $\Lambda(\theta)+\theta\delta_1=-\sqrt{\delta}\Lambda'(\gamma)+(\gamma-\sqrt{\delta})\delta_1+O(\delta)=-\sqrt{\delta}\Lambda'(\gamma)+O(\delta)<0$ for $\delta$ small enough, thus the upper bound in \eqref{subeq: upper bound 3} holds. Second, we see that
\begin{multline*}
\sup_{\eta\in\mathbb{R}}\left(\eta\left(\frac{1}{\gamma}-\delta_2\right)-\left(\frac{1}{\Lambda'(\gamma)\gamma}+\epsilon\right)(\Lambda(\eta)+\eta\delta_1)\right)\\
\geq (\gamma-\sqrt{\delta})\left(\frac{1}{\gamma}-\delta_2\right)-\left(\frac{1}{\Lambda'(\gamma)\gamma}+\epsilon\right)(\Lambda(\gamma-\sqrt{\delta})+(\gamma-\sqrt{\delta})\delta_1).
\end{multline*}
So, we can conclude that 
$$
(\gamma-\sqrt{\delta})\left(\frac{1}{\gamma}-\delta_2\right)-\left(\frac{1}{\Lambda'(\gamma)\gamma}+\delta^{1/4}\right)(\Lambda(\gamma-\sqrt{\delta})+(\gamma-\sqrt{\delta})\delta_1)=1+\delta^{3/4}\Lambda'(\gamma)+O(\delta)>1
$$
for $\delta$ small enough, thus the expression in \eqref{subeq: upper bound 3} converges to 0 as $N\to\infty$.
\end{proof}
\begin{proof}[Proof of Theorem \ref{thm: convergence maximum waiting time}]
First, to prove a lower bound, we see that
$$
\max_{i\leq N}W_i(\infty)\geq_{st.} \max_{i\leq N}\sum_{j=1}^{\big\lfloor \frac{1}{(\Lambda'(\gamma)\gamma)}\log N\big\rfloor}(S_i(j)-A(j)).
$$
Thus, combining this inequality with the result from Lemma \ref{lem: lower bound}, we see that 
$$
\liminf_{N\to\infty}\probability*{\max_{i\leq N}W_i(\infty)>\frac{1}{\gamma}\log N+x\sqrt{\log N}}\geq \probability*{\sigma_A\sqrt{\frac{1}{\Lambda'(\gamma)\gamma}}X>x}.
$$
Second, by using the union bound of the types as given in \eqref{eq: upper bound} and explained in Section \ref{sec: heuristic analysis}, we get from Lemmas \ref{lem: supremum 0 t}, \ref{lem: t-e t+e}, and \ref{lem: t+e infty}, with  $t_1^{(N)}=\left(\frac{1}{\Lambda'(\gamma)\gamma}-\epsilon\right)\log N$ and $t_3^{(N)}=\left(\frac{1}{\Lambda'(\gamma)\gamma}+\epsilon\right)\log N$, that
\begin{align*}
    &\limsup_{N\to\infty}\probability*{\max_{i\leq N}W_i(\infty)>\frac{1}{\gamma}\log N+x\sqrt{\log N}}\\
    &\quad\leq \limsup_{N\to\infty}\probability*{\max_{i\leq N}\sup_{t_1^{(N)}\leq k<t_3^{(N)}}\sum_{j=1}^k(S_i(j)-A(j))>\frac{1}{\gamma}\log N+x\sqrt{\log N}}\\
&\quad\leq  \probability*{\sigma_A\sqrt{\frac{1}{\Lambda'(\gamma)\gamma}-\epsilon}X_1+\sigma_A\sqrt{2\epsilon}\left|X_2\right|>x}.
\end{align*}
Finally, we have that 
$$
\probability*{\sigma_A\sqrt{\frac{1}{\Lambda'(\gamma)\gamma}-\epsilon}X_1+\sigma_A\sqrt{2\epsilon}\left|X_2\right|>x}\overset{\epsilon\downarrow 0}{\longrightarrow} \probability*{\sigma_A\sqrt{\frac{1}{\Lambda'(\gamma)\gamma}}X>x}.
$$
\end{proof}


\begin{thebibliography}{10}

\bibitem{asmussen2003applied}
S{\o}ren Asmussen.
\newblock {\em Applied probability and queues}, volume~2.
\newblock Springer, 2003.

\bibitem{baccelli1985two}
Fran{\c{c}}ois Baccelli.
\newblock {Two parallel queues created by arrivals with two demands: The M/G/2
  symmetrical case}.
\newblock {\em Technical report RR--0426, INRIA}, 1985.

\bibitem{baccelli1989queueing}
Fran{\c{c}}ois Baccelli and Armand~M. Makowski.
\newblock Queueing models for systems with synchronization constraints.
\newblock {\em Proceedings of the IEEE}, 77(1):138--161, 1989.

\bibitem{dembo2009large}
Amir Dembo and Ofer Zeitouni.
\newblock {\em Large deviations techniques and applications}, volume~38.
\newblock Springer Science \& Business Media, 2009.

\bibitem{donsker1951invariance}
Monroe~David Donsker.
\newblock {\em An invariance principle for certain probability limit theorems},
  volume~6.
\newblock Memoirs of the American Mathematical Society, 1951.

\bibitem{flatto1984two}
Leopold Flatto and Sann Hahn.
\newblock {Two parallel queues created by arrivals with two demands I}.
\newblock {\em SIAM Journal on Applied Mathematics}, 44(5):1041--1053, 1984.

\bibitem{de2007extreme}
Laurens {\VAN{Haan}{de}{de}}~Haan and Ana Ferreira.
\newblock {\em Extreme value theory: an introduction}.
\newblock Springer Science \& Business Media, 2006.

\bibitem{haji1971relation}
Rasoul Haji and Gordon~F Newell.
\newblock A relation between stationary queue and waiting time distributions.
\newblock {\em Journal of Applied Probability}, 8(3):617--620, 1971.

\bibitem{de1988fredholm}
Stephanus J.~de Klein.
\newblock {\em Fredholm integral equations in queueing analysis}.
\newblock PhD thesis, Rijksuniversiteit Utrecht, 1988.

\bibitem{ko2004response}
Sung-Seok Ko and Richard~F. Serfozo.
\newblock {Response times in M/M/s fork-join networks}.
\newblock {\em Advances in Applied Probability}, 36(3):854--871, 2004.

\bibitem{lindley1952theory}
David~V Lindley.
\newblock The theory of queues with a single server.
\newblock In {\em Mathematical Proceedings of the Cambridge Philosophical
  Society}, volume~48, pages 277--289. Cambridge University Press, 1952.

\bibitem{MeijerSchol2021}
Mirjam Meijer, Dennis Schol, Willem {\VAN{Jaarsveld}{van}{van}}~Jaarsveld,
  Maria Vlasiou, and Bert Zwart.
\newblock {Extreme-value theory for large fork-join queues, with applications
  to high-tech supply chains}.
\newblock {\em https://arxiv.org/abs/2105.09189}, 2021.

\bibitem{michel1981constant}
Reinhard Michel.
\newblock {On the constant in the nonuniform version of the Berry-Ess\'een
  theorem}.
\newblock {\em Zeitschrift f{\"u}r Wahrscheinlichkeitstheorie und verwandte
  Gebiete}, 55(1):109--117, 1981.

\bibitem{nelson1988approximate}
Randolph Nelson and Asser~N Tantawi.
\newblock Approximate analysis of fork/join synchronization in parallel queues.
\newblock {\em IEEE Transactions on Computers}, 37(6):739--743, 1988.

\bibitem{scholMOR}
Dennis Schol, Maria Vlasiou, and Bert Zwart.
\newblock Large fork-join queues with nearly deterministic arrival and service
  times.
\newblock {\em Mathematics of Operations Research}, 47(2):1335--1364, 2021.

\bibitem{schol2022tail}
Dennis Schol, Maria Vlasiou, and Bert Zwart.
\newblock Tail asymptotics for the delay in a brownian fork-join queue.
\newblock {\em arXiv preprint arXiv:2208.04796}, 2022.

\bibitem{wright1992two}
Paul~E. Wright.
\newblock Two parallel processors with coupled inputs.
\newblock {\em Advances in Applied Probability}, 24(4):986--1007, 1992.

\end{thebibliography}
\end{document}